\documentclass[11pt]{article}
\usepackage{amsfonts}
\usepackage{latexsym, amssymb, amsmath, amscd, amsfonts, epsfig, graphicx, colordvi}
\usepackage{ifpdf}
\usepackage{graphicx}
\newtheorem{thm}{Theorem}[section]

\def\qed{\nopagebreak\hfill{\rule{4pt}{7pt}}
\medbreak}

\def\qed{\nopagebreak\hfill{\rule{4pt}{7pt}}
\medbreak}

\makeatletter
\def\ExtendSymbol#1#2#3#4#5{\ext@arrow 0099{\arrowfill@#1#2#3}{#4}{#5}}
\makeatother

\title{ Involutions for Rogers-Ramanujan-Gordon Type Identities with Parity Restrictions}
\author{William Y.C. Chen\raisebox{5pt}{\scriptsize 1},
Doris D. M. Sang\raisebox{5pt}{\scriptsize 2}, and Diane Y. H. Shi\raisebox{5pt}{\scriptsize 3}}
\date{Center for Combinatorics, LPMC-TJKLC\\
 Nankai University\\
Tianjin 300071, P.R. China \\
\vspace{15pt}
\raisebox{5pt}{\scriptsize 1\,}chen@nankai.edu.cn,
\raisebox{5pt}{\scriptsize 2\,}sdm@cfc.nankai.edu.cn, \raisebox{5pt}{\scriptsize 3\,}yahuishi@gmail.com}

\begin{document}
\maketitle
\noindent {\bf Abstract.} We find involutions for
three Rogers-Ramanujan-Gordon type
identities obtained by Andrews on the generating functions
for partitions with part difference and parity restrictions.



\section{ Introduction}

This paper is concerned with combinatorial interpretations of three Rogers-Ramanujan-Gordon type
identities obtained  by Andrews  on the generating functions of
partitions with parity restrictions.
 Recall that Gordon \cite{gor61} found a combinatorial generalization of the Rogers-Ramanujan identities \cite{and76,rog1894}, which has been called the
 Rogers-Ramanujan-Gordon identity, see Andrews \cite{and66}.

\begin{thm}\label{thm1.2} Let $B_{k,a}(n)$ denote the number of partitions
of $n$ for the form $b_1 + b_2 + \cdots + b_j$, where $b_i \geq b_i+1$,
$b_i-b_{i+k-1}\geq2$ and at most $a-1$ of the $b_i$ are equal to $1$
and $1\leq a\leq k$. Let $A_{k,a}(n)$ denote the number of
partitions of $n$ into parts $\not \equiv0,\pm a \;(\mbox{mod}\ 2k + 1)$.
Then for all $n \geq0$,
 \begin{equation}\label{ab}
 A_{k,a}(n) = B_{k,a}(n).
 \end{equation}
\end{thm}

It is easy to derive the generating function for $A_{k,a}(n)$,
\begin{equation}\label{af}
\sum_{n=0}^{\infty}A_{k,a}(n)q^n
=\frac{(q^{2k+1};q^{2k+1})_{\infty}
(q^a;q^{2k+1})_{\infty}(q^{2k+1-a};q^{2k+1})_{\infty}}{(q;q)_{\infty}}.
\end{equation}
In order to prove (\ref{ab}), Gordon has shown that  $B_{k,a}$ has
the same generating function  as $A_{k,a}$, that is,
\begin{equation}\label{eqbf}
\sum_{n=0}^{\infty}B_{k,a}(n)q^n
=\frac{(q^{2k+1};q^{2k+1})_{\infty}
(q^a;q^{2k+1})_{\infty}(q^{2k+1-a};q^{2k+1})_{\infty}}{(q;q)_{\infty}}.
\end{equation}
In fact, Gordon established the above identity by giving an involution
for the following reformulation of \eqref{eqbf} by
multiplying both sides by $(q;q)_\infty$ and applying the Jacobi triple
product identity to the right hand side,
\begin{equation}\label{ebf}
(q;q)_{\infty}\sum_{n=0}^{\infty}B_{k,a}(n)q^n
=\sum_{n=-\infty}^{\infty}(-1)^nq^{(k+\frac{1}{2})n^2+(k-a+\frac{1}{2})n}.
\end{equation}
We shall give a description of
Gordon's involution in Section 2.
Andrews \cite{and74} found  a generating function proof of
the Rogers-Ramanujan-Gordon identity \eqref{eqbf}
in the following equivalent form:
\[\sum_{n_1,\ldots,n_{k-1}\geq
0}\frac{q^{N_1^2+N_2^2+\cdots+N_{k-1}^2+N_a+N_{a+1}+\cdots+N_{k-1}}}{(q;q)_{n_1}(q;q)_{n_2}\cdots(q;q)_{n_{k-1}}}
=\frac{(q^{2k+1};q^{2k+1})_{\infty}(q^a;q^{2k+1})_{\infty}(q^{2k+1-a};q^{2k+1})_{\infty}}{(q;q)_{\infty}},\] where
\[N_j=n_j+n_{j+1}+\cdots+n_{k-1},\]
and \[(a;q)_n=(a)_n=(1-a)(1-aq)\cdots(1-aq^{n-1}).\]

Andrews \cite{and07} further considered the parity restrictions on partitions in connection with the Rogers-Ramanujan-Gordon identity. He has derived three
identities in the spirit of the relation \eqref{eqbf}, which we call Rogers-Ramanujan-Gordon type identities.  The three identities depend on the
parities of $k$ and $a$.

\begin{thm}\label{thm1.3}
Suppose $k\geq a\geq1$ are integers such that $k$ and $a$ are both even. Let
$W_{k,a}(n)$ denote the number of partitions enumerated by
$B_{k,a}(n)$ with further restriction that even parts appear an
even number of times. Then we have
\begin{equation}\label{equ1.1}\sum_{n\geq 0}W_{k,a}(n)q^n=\frac{(-q;q^2)_{\infty}(q^a;q^{2k+2})_{\infty}
(q^{2k+2-a};q^{2k+2})_{\infty}(q^{2k+2};q^{2k+2})_{\infty}}{(q^2;q^2)_{\infty}}.\end{equation}
\end{thm}

\begin{thm}\label{thm1.4}
Suppose $k\geq a\geq1$ are integers such that $k$ and $a$ are both odd. Let
$W_{k,a}(n)$ denote the number of those partitions enumerated by
$B_{k,a}(n)$ with further restriction that even parts appear an
even number of times. Then for all $n\geq0$, we have
\begin{equation}\label{equ1.2}\sum_{n\geq 0}W_{k,a}(n)q^n=\frac{(q^2;q^4)_{\infty}(q^a;q^{2k+2})_{\infty}
(q^{2k+2-a};q^{2k+2})_{\infty}(q^{2k+2};q^{2k+2})_{\infty}}{(q;q)_{\infty}}.\end{equation}
\end{thm}

\begin{thm}\label{thm1.5}
Suppose $k\geq a\geq1$ and  $k$ odd and $a$ even. Let $\overline{W}_{k,a}(n)$ denote the
number of those partitions enumerated by $B_{k,a}(n)$ with further
restriction that odd parts appear an even number of times. Then for
all $n\geq0$, we have
\begin{equation}\label{equ1.3}\sum_{n\geq0}\overline{W}_{k,a}(n)q^n=\frac{(q^a;q^{2k+2})_{\infty}
(q^{2k+2-a};q^{2k+2})_{\infty}(q^{2k+2};q^{2k+2})_{\infty}}{(-q;q^2)_{\infty}(q;q)_{\infty}}.\end{equation}
\end{thm}

Andrews \cite{and07} gave algebraic proofs of the
above three identities.
In this paper, we shall give three involutions  for   Andrews' identities.
Gordon's involution will play a role in our involutions for the purpose of partial cancelations.

\section{Gordon's involution}

In this section we give an overview of Gordon's involution
for the Rogers-Ramanujan-Gordon identity. Recall that
 $B_{k,a}(N)$ denotes the number of partitions
of $N$ for the form $b_1 + b_2 + \cdots +b_m $, where $b_i \geq b_i+1$,
$b_i-b_{i+k-1}\geq2$ and at most $a-1$ of the $b_i$ are equal to $1$
and $1\leq a\leq k$. Gordon's involution is concerned with the following
relation
\begin{equation}\label{geq}(q;q)_{\infty}\sum_{N=0}^{\infty}B_{k,a}(N)q^N
=\sum_{n=-\infty}^{\infty}(-1)^nq^{(k+1/2)n^2+(k-a+1/2)n}.\end{equation}

Let us consider the set $\mathcal{P}_{k,a}$ of pairs of partitions $(A|B)$, where $A$ is a signed partition with distinct parts and $B$ is a partition enumerated by $B_{k,a}$. To be precise, each part in the partition $A$ carries a minis sign, and $B$ is a partition of the form $N=b_1+b_2+\cdots+b_m$, where $b_i\geq b_{i+1}$, $b_i\geq b_{i+k-1}+2$ and  $b_{m-a+1}\geq2$. Clearly, the set $\mathcal{P}_{k,a}$ corresponds to the
left hand side of (\ref{geq}). Gordon's involution is based on the set $\mathcal{P}_{k,a}$.
As will be shown, the generating function of the fixed points of Gordon's involution equals
\[\sum_{n=-\infty}^{\infty}(-1)^nq^{(k+\frac{1}{2})n^2+(k-a+\frac{1}{2})n},\]
which is the right hand side of \eqref{geq}.

Gordon's involution consists of the following steps.

\noindent Step 1.
We compare the largest parts of $A$ and $B$, if $a_1> b_1$  then set \[(A'|B')= (b_1,a_1,\ldots,a_l|b_2,\ldots,b_m).\]
 In this case, the resulting pair $(A'|B')$ is still in $\mathcal{P}_{k,a}$.
  In the case $a_1\leq b_1$, if there exists $1\leq i \leq k$ such that
\begin{equation}\label{equ2.3}
a_1=b_1=\cdots=b_{i-1}=b_i+1=\cdots=b_{k-1}+1,\end{equation}
then go to Step 2. Otherwise, let
\[(A'|B')=(a_2,\ldots,a_l|a_1,b_1,b_2,\ldots,b_m).\]
It can be checked that the resulting pair $(A'|B')$ is in $P_{k,a}$

\noindent Step 2. Let $U_i$ denote the  set of  partition pairs $(A|B)$
for which the condition
\eqref{equ2.3} holds, and let $U$ denote the union of $U_i$.   We shall further classify the sets $U_i$.
Let us define the three  numbers $p$, $q$ and  $r$ as
follows:

(1) $p=a_l$.

(2) $q$ is the largest integer such that
\[a_1-a_2=a_2-a_3=\cdots=a_{q-1}-a_q=1.\]

(3)  $r$ is the
largest integer such that
\[b_{k-1}-b_{2(k-1)}=b_{2(k-1)}-b_{3(k-1)}=\cdots=b_{(r-1)(k-1)}-b_{r(k-1)}=2.\]

 For $U_2,\ldots,U_{k-1}$, there is a fourth
parameter number $s$, which is  defined to be the largest integer such
that
\[b_{i-1}-b_{i-1+(k-1)}=b_{i-1+(k-1)}-b_{i-1+2(k-1)}=\cdots=b_{i-1+(s-2)(k-1)}-b_{i-1+(s-1)(k-1)}=2.\]

Now we divide each of $U_1,U_k$ into three classes, and each of $U_2,\ldots,U_{k-1}$
into four classes. For ${U}_1,{U}_k$, we consider the minimal of the three numbers $p,q$ and $r$, and let $n=min\{p,q,r\}$. For $2\leq i\leq k-1$, we consider the minimal of the four numbers $p$, $q$, $r$ and $s$, and let $n=min\{p,q,r,s\}$. We now divide the sets $U_i$ into classes according to
the conditions on the parameters $p$, $q$, $r$ and $s$:

\[{U}^1_1:\ \ \ p=n,\ q\geq n, r\geq n,\]
\[{U}^2_1:\ \ \ p>n,\ q=n,\ r\geq n,\]
\[{U}^3_1:\ \ \ p>n,\ q>n,\ r=n,\]
\[{U}^1_k:\ \ \ p=n,\ q\geq n,\ r\geq n,\]
\[{U}^2_k:\ \ \ p>n,\ q\geq n,\ r=n,\]
\[{U}_k^4:\ \ \ p>n,\ q=n,\ r>n,\]
\[{U}^1_i:\ \ \ p=n,\ q\geq n,\ r\geq n,\]
\[{U}^2_i:\ \ \ p>n,\ q\geq n,\ r\geq n, s=n\]
\[{U}^3_i:\ \ \ p>n,\ q>n,\ r=n,\ s>n,\]
\[{U}_i^4:\ \ \ p>n,\ q=n,\ r\geq n,\ s>n.\]

\noindent Step 3. We now define three maps $
\alpha$, $\beta$, $\gamma$ and their inverses.
Let  \[\pi=(A|B)=(a_1,\ldots,a_l|b_1,\ldots,b_m)\in U^2_1,\]
define
\[\alpha(\pi)=(a_1-1,a_2-1,\ldots,a_n-1,a_{n+1},\ldots,a_l,n|b_1,b_2,\ldots,b_m),\]
It is clear that $\alpha(\pi)$ belongs  $U^1_k$, unless $l=n$,
$a_n-1=n$, and $\pi$ of the following form
\[(2n,2n-1,\ldots,n+1|(2n-1)^{k-1}\ (2n-3)^{k-1}\ldots\ 1^{k-1}).\]
In this case, the weight of $\pi$ equals $(k+\frac{1}{2})n^2+\frac{1}{2}n$.

Conversely, if $\pi\in U^1_k$,  we have
\[\alpha^{-1}(\pi)=(a_1+1,\ldots,a_n+1,a_{n+1},\ldots,a_{l-1}|b_1,\ldots,b_m),\]
which belongs  $U^2_1$, unless $l=n$. For the case $l=n$  we have
\[\pi=(2n-1,2n-2,\ldots,n|(2n-1)^{k-1}\ (2n-3)^{k-1}\ldots\ 1^{k-1}),\]
which has weight $(k+\frac{1}{2})n^2-\frac{1}{2}$. In this case,
it can be seen that the part $1$  appears $k-1$ times. On the
other hand, by the condition in the theorem, the element $1$ may appear at most $a-1$ times and $a\leq k$, and so the two exceptional cases can happen only when $a=k$.

If $\pi\in U^1_{i}$ for some $1\leq i\leq k-1$,  we put
\begin{align*}\beta(\pi)=(a_1,\dots,a_{l-1}|b_1,\ldots,b_{i-1},b_i+1,b_{i+1},\ldots,b_{i-1+(k-1)},\\
b_{i+(k-1)}+1,\ldots,b_{i+(n-1)(k-1)}+1,\ldots).\end{align*}
Clearly, $\beta(\pi)\in U^2_{i+1}$.

Conversely,  if $\pi\in U^2_{i}$ for some $2\leq i\leq k$, then we have
\begin{equation}\label{bpi}
\beta^{-1}(\pi)=
(a_1,\ldots,a_l,n|b_1,\ldots,b_{i-1},b_i-1,\ldots,b_{i+(n-1)(k-1)}-1,\ldots)
\end{equation}
belongs $U^1_{i-1}$ unless $i=k-a$. In the case of $i=k-a$, the element $1$ appears more than $a-1$ times in the  partition in (\ref{bpi}).
We deduce that
\[\pi=(2n,2n-1,\ldots,n+1|(2n)^{k-a}\ (2n-1)^{a-1}\ (2n-2)^{k-a}\ldots),\] which
 has weight
$(k+\frac{1}{2})n^2+(k-a+\frac{1}{2})n.$

If $\pi\in U^4_i$ for some $2\leq i\leq k$, we define
\begin{align*}\gamma(\pi)=(b_1,a_1-1,&\ldots,a_n-1,a_{n+1},\ldots,a_l|\\
&b_2,\ldots,b_{k-1},b_k+1,b_{k+1},\ldots,b_{2k-1}+1,\ldots,b_{k+(n-1)(k-1)}+1,\ldots).\end{align*}
Then $\gamma(\pi)\in U^3_{i-1}$.

Conversely if $\pi\in U^3_{i}$, for some $1\leq i\leq k-1$,
then
\begin{equation}\label{gpi}
\gamma^{-1}(\pi)=
(a_2+1,\ldots,a_{n+1}+1,\ldots,a_l
|a_1,b_1,\ldots,b_{k-1}-1,\ldots,b_{(n-1)(k-1)}-1,\ldots),\end{equation}
 belongs $U^4_i$
  unless $a=i$.  When $a=i$, the element $1$ appears more than $a-1$ times in the partition in (\ref{gpi}). So we are led to the fixed point
\[\pi=(2n-1,2n-2,\ldots
,n|(2n-1)^{a-1}\ (2n-2)^{k-a}\ (2n-3)^{a-1}\ldots),\] which has weight $(k+\frac{1}{2})n^2-(k-a+\frac{1}{2})n$.

In summary, we get two types of fixed points
\[(2n,2n-1,\ldots,n+1|(2n)^{k-a}\ (2n-1)^{a-1}\ (2n-2)^{k-a}\ldots)\] and
\[(2n-1,2n-2,\ldots
,n|(2n-1)^{a-1}\ (2n-2)^{k-a}\ (2n-3)^{a-1}\ldots),\]
where $n\geq 0$.

For example, let $k=a=3$ and $N=17$. Then the correspondence that is
not covered by Step 1   consists of the following maps.

(1)  $\alpha \colon U_1^2 \rightarrow U_3^1$:
\begin{center}
\begin{tabular}{llrl}
 & $(6\ |\ 5,5,1)$ & $\begin{CD} @>\alpha>> \end{CD}$ &$(5\ |\ 5,4,3)$\\[3pt]
&$(5,3\ |\ 4,4,1)$& &$(4,3,1\ |\ 4,4,1)$\\[3pt]
&$(5,2\ |\ 4,4,2)$& &$(4,2,1\ |\ 4,4,2)$\\[3pt]
&$(5,2\ |\ 4,4,1,1)$ && $(4,2,1\ |\ 4,4,1,1)$\\[3pt]
&$(5\ |\ 4,4,2,2)$ && $(4,1\ |\ 4,4,2,2)$.
\end{tabular}
\end{center}

(2)  $\beta\colon U_1^1 \rightarrow U_2^2$:
\begin{center}
\begin{tabular}{llrl}
 & $(6,1\ |\ 5,5)$ & $\begin{CD} @>\beta>> \end{CD}$ &$(6\ |\ 6,5)$\\[3pt]
&$(5,3,1\ |\ 4,4)$& &$(5,3\ |\ 5,4)$\\[3pt]
&$(5,2,1\ |\ 4,4,1)$& &$(5,2\ |\ 5,4,1)$\\[3pt]
&$(5,1\ |\ 4,4,2,1)$ && $(5\ |\ 5,4,2,1)$\\[3pt]
&$(4,3,2,1\ |\ 3,3,1)$ && $(4,3,2\ |\ 4,3,1)$.
\end{tabular}
\end{center}

(3) $\beta \colon U_2^1 \rightarrow U_3^2$:
\begin{center}
\begin{tabular}{llrl}
 & $(5,2,1\ |\ 5,4)$ & $\begin{CD} @>\beta>> \end{CD}$ &$(5,2\ |\ 5,5)$\\[3pt]
&$(5,1\ |\ 5,4,2)$& &$(5\ |\ 5,5,2)$\\[3pt]
&$(5,2\ |\ 4,4,2)$& &$(5\ |\ 5,5,1,1)$\\[3pt]
&$(4,3,2,1\ |\ 4,3)$ && $(4,3,2\ |\ 4,4)$\\[3pt]
&$(4,3,1\ |\ 4,3,2)$ && $(4,3\ |\ 4,4,2)$\\[3pt]
&$(4,3,1\ |\ 4,3,1,1)$&&$(4,3\ |\ 4,4,1,1)$\\[3pt]
&$(4,2,1\ |\ 4,3,2,1)$&&$(4,2\ |\ 4,4,2,1).$
\end{tabular}
\end{center}

(4) $\gamma\colon U_2^4 \rightarrow U_1^3$:
\begin{center}
\begin{tabular}{llrl}
 & $(5\ |\ 5,4,3)$ & $\begin{CD} @>\gamma>> \end{CD}$ &$(5,4\ |\ 4,4).$
\end{tabular}
\end{center}
In this example, there are no fixed points.

\section{Combinatorial proof of Theorem \ref{thm1.3} }

In this section, we give an involution for Andrews' Theorem \ref{thm1.3}.
By the relation \[\frac{(-q;q^2)_{\infty}}{(q^2;q^2)_{\infty}}=\frac{(q^2;q^4)_{\infty}}{(q;q)_{\infty}},\] and the Jacobi triple product identity,  Theorem \ref{thm1.3} can be restated
as
\begin{equation} \label{a}
(q;q)_{\infty}\sum_{n=0}^{\infty}W_{k,a}(n)q^n=(q^2;q^4)_{\infty}
\sum_{n=-\infty}^{\infty}(-1)^nq^{(k+1)n^2+(k+1-a)n}.\end{equation}
 As the first step, we shall transform the partition pairs $(A|B)$
 corresponding  the left hand side of identity \eqref{a} into  triples of partitions $(A'|B'|E)$  such that partition pairs $(A'|B')$  satisfy the conditions for Gordon's involution for some $k$ and $a$. Then  we apply Gordon's involution on $(A'|B')$ and keep $E$ unchanged. It will be shown that the right hand side of (\ref{a})
 equals the generating function for the fixed points $(A'|B'|E)$ of this involution.

\noindent
{\it Combinatorial Proof of  Theorem \ref{thm1.3}.}
The left hand side of (\ref{a})
 can be interpreted as  the generating function of
 pairs of partitions $(A|B)$, where
$A= (a_1,a_2,\ldots)$ is a signed  partition with distinct parts and $B= (b_1,b_2,\ldots)$ is a partition enumerated by  $W_{k,a}$. We denote the set of such partition pairs by $\mathcal{G}_{k,a}$

We shall  construct an involution on $\mathcal{G}_{k,a}$ which leads to a combinatorial
proof of (\ref{a}). It consists of five steps. The objective of the first two steps
is transform a pair
of partitions $(A|B)$ in $\mathcal{G}_{k,a}$ into a triple of partitions $(A'|B'|E)$ such that
$E$ is a signed partition with parts congruent to $2$ modulo $4$ and $B'$ is a
partition with odd parts distinct.  Step 3 gives an involution to cancel out certain
partition pairs so that we can apply Gordon's involution to the remaining partition
pairs, which is Step 4. In Step 5, we still need an involution to reach the fixed points.

\noindent Step 1. Let $(A|B)$ be a partition pair in $\mathcal{G}_{k,a}$ as given above.
 If $b_i$ is an odd part in $B$ and there is a part
in $A$ that equals $b_i$, then we remove $b_i$ from $B$ and the equal part from $A$  to form a part in $E$ that is congruent to $2$ modulo $4$.
Since the part in $A$ has a minus sign, the resulting part in $E$ also has
a minus sign. Repeating this procedure until there are no common odd parts
 in $A$ and $B$.

\noindent Step 2.
Combining two equal parts into one part in $B'$ until we get a partition
 with odd parts distinct. For the sake of notational convenience, we still use $B'$ to
denote the resulting partition.
We denote the set of the
 triples of partitions $(A'|B'|E)$ by $\mathcal{G}_{k,a}'$.

\noindent Step 3. We compare the largest parts in $A'$ and $B'$, namely,
$a'_1$ and $b'_1$. If $b'_1> a'_1$, then we move  $b'_1$ to
$A'$. This operation leads to a triple of partitions in $\mathcal{G}_{k,a}'$.
If $b'_1\leq a'_1$, we move $a'_1$ to  $B'$, except for the following case. To describe the exceptional case,
let us rewrite $B'$ as a pair of partitions $(C|D)$, where
$C=(c_1,c_2,\ldots)$ consists of the odd parts of $B'$ and $D=(d_1,d_2,\ldots)$ consists of the even parts
of $B'$. The conditions are as follows

(1) $a_1=d_1=d_2=\cdots=d_{i-1}=d_i+2=\cdots=d_{\frac{k-2}{2}}+2$, or

(2) $a_1=d_1=d_2=\cdots=d_{i-1}=d_i+2=\cdots=d_{\frac{k-4}{2}}+2$,
and there is an odd part $c_j$ such that  $c_j=a_1/2$ or $ a_1/2 -1$.

It can be checked that if either (1) or (2) holds, then the partition triples
obtained from $(A'|B'|E)$ by moving $a'_1$ to $B$ no longer belongs
$\mathcal{G}_{k,a}'$. Moreover, one sees that if neither (1) nor (2) holds, then
the triple of partitions obtained by the above operation is still in $\mathcal{G}_{k,a}'$.

\noindent Step 4. We shall give an involution for triples of partitions in $\mathcal{G}_{k,a}'$
that satisfy the above conditions (1) or (2). Note that the odd parts in $A'$ and $C$ are all
distinct. Let us compare the largest odd parts in $A'$ and $C$.
When the larger one is in $A'$, we  move it to $C$. Similarly, when the
larger one is in $C$, we move
the largest part of $C$ to $A'$. To be precise, if only one of $A'$ and $C$
 contains odd parts, then we move the largest odd part to the other partition.
 So there is only one case that it is impossible to make any move, that is,
 neither $A'$ nor $C$ contains odd parts.

Till now, we have left the triples of partitions in $\mathcal{G}_{k,a}'$
such that $(A'|B')$ is a pair of partitions with even parts.

\noindent Step 5. For remaining partition pairs $(A'|B')$, we apply Gordon's involution  to partition pairs $(A'|B')$ with $k$ and $a$ replaced by $\frac{k}{2}$ and $\frac{a}{2}$ respectively and with every part doubled.
Then we get two types of  fixed points
\begin{equation}\label{1t}
(A'|B'|E)=(4n,4n-2,\ldots,2n+2|(4n)^{\frac{k-a}{2}},
(4n-2)^{\frac{a-2}{2}},(4n-4)^{\frac{k-a}{2}}\ldots|E),
\end{equation}
and
\begin{equation}\label{2t}
(A'|B'|E)=(4n-2,\ldots,2n|
(4n-2)^{\frac{a-2}{2}},(4n-4)^{\frac{k-a}{2}},
(4n-6)^{\frac{a-2}{2}}\ldots|E).
\end{equation}
For the partition triples in (\ref{1t}), we have
\[|A'|+|B'|=(k+1)n^2+(k-a+1)n,\]
whereas for the partition triples in (\ref{2t}), we have
\[|A'|+|B'|
=(k+1)n^2-(k-a+1)n.\]
Since  $E$ is any signed partition with distinct parts congruent to $2$
modulo $4$, we see that the generating function of the fixed points $(A'|B'|E)$ equals
\[(q^2;q^4)_{\infty}
\sum_{n=-\infty}^{\infty}(-1)^nq^{(k+1)n^2+(k+1-a)n}.\]
This completes the combinatorial proof of Theorem \ref{thm1.3}. \qed

For example, let $k=a=6$ and $n=44$, and let $(A|B)=(10,8,5\ |\ 5,4^4)$. The first two steps are as follows:
\[(10,8,5\ |\ 5,4^4)\Rightarrow(10,8\ |\ 4^4\ |\ 5,5)\Rightarrow(10,8\ |\ 8,8\ |\ 10).\]
 Now we see that the partition $(10,8\ |\ 8,8\ |\ 10)$   belongs  $U^3_1$
 as defined in Section 2 if every part for Gordon's involution is doubled.
   Using the map $\gamma^{-1}$ also with every part doubled,  we obtain
    the partition $(10\ |\ 10,8,6\ |\ 10)$  from the partition $(10,5\ |\ 5,5,5,4,4,3,3)$ by the following steps
\[(10,5\ |\ 5,5,5,4,4,3,3)\Rightarrow(10\ |\ 5,5,4,4,3,3\ |\ 5,5)\Rightarrow(10\ |\ 10,8,6\ |\ 10).\]

\section{Combinatorial proof of  Theorem \ref{thm1.4} }

In this section, we give a combinatorial proof of Theorem \ref{thm1.4} which
also involves Gordon's involution. Using the relation \[\frac{(q^2;q^4)_{\infty}}{(q;q)_{\infty}}=\frac{(-q;q^2)_{\infty}}{(q^2;q^2)_{\infty}},\] and the Jacobi triple product identity,
Theorem \ref{thm1.4} can be rewritten as follows
\begin{equation}\label{b}
(q^2;q^2)_{\infty}\sum_{n=0}^{\infty}W_{k,a}(n)q^n=(-q;q^2)_{\infty}
\sum_{n=-\infty}^{\infty}(-1)^nq^{(k+1)n^2+(k+1-a)n}.\end{equation}
Suppose that $(A|B)$ is a partition pair corresponding to the left hand side of the identity \eqref{b}. We first transform it into a triple of partitions $(A|C'|D')$. Then we
establish an involution on the set of these triples such that the fixed
points  have the generating function equals the right hand side of \eqref{b}.

\noindent {\it Combinatorial proof of the Theorem\ref{thm1.4}}.
The left hand side
of (\ref{b}) can be interpreted
 as   the generating function of pairs of partitions $(A|B)$, where $A=(
a_1,a_2,\ldots)$ is a
signed partition with distinct even parts, and $B=( b_1,b_2,\ldots)$ is a partition enumerated by $W_{k,a}$.  We denote the set of such
pairs $(A|B)$ by $\mathcal{H}_{k,a}$.

We proceed to construct an
involution on the set $\mathcal{H}_{k,a}$.
 It consists of five steps. The objective of the first step
is transform a pair
of partitions $(A|B)$ in $\mathcal{H}_{k,a}$ into a triple of partitions $(A|C|D)$ such that
$C$ is a  partition with even parts and $D$ is a
partition with distinct odd parts.  Step 2 gives an involution to cancel out certain
partition pairs. In Step 3, we split the parts in $C$ that are congruent to $2$ modulo $4$ but not a double of some part in $D$ into two equal parts and
move them to $D$. In Step 4, we  apply   Gordon's involution. Step 5 gives an operation on partitions that enables us to compute the generating function of the fixed points.

\noindent Step 1. Combine two equal  parts in $B$ to form an even part. Repeating this
operation to generate a partition $C=(c_1, c_2, \ldots)$  with even parts. Let $D=(d_1, d_2, \ldots)$ be the
 partitions consisting of the remaining parts in $B$. Since the even parts in $B$ appear an even number of times, $D$ is a partition with odd parts.
 We denote  the set of the  triples  $(A|C|D)$  by $\mathcal{H}_{k,a}'$.

\noindent Step 2. Now we compare the largest parts in $A$ and $C$.  If $c_1> a_1$, then we move the $c_1$ to $A$. This gives a triple of partitions  in $\mathcal{H}_{k,a}'$. If $c_1\leq a_1$, we move $a_1$ to  $C$ to form a triple  in  $\mathcal{H}_{k,a}'$ except for the following cases:

(1) $a_1=c_1=c_2=\cdots=c_{i-1}=c_i+2=\cdots=c_{\frac{k-1}{2}}+2$, or

(2) $a_1=c_1=c_2=\cdots=c_{i-1}=c_i+2=\cdots=c_{\frac{k-3}{2}}+2$,
and there is an odd part $d_j=a_1/2$ or $(a_1-2)/2$.

\noindent Step 3. For the triples of partitions $(A|C|D)$
in $\mathcal{H}_{k,a}'$ that satisfy the above
exceptional conditions, we  construct an involution to further cancel out certain
triples. First, we need to move some parts of $C$ to $D$.
Consider the parts in $C$ that are congruent to $2$ modulo $4$. If
a part $4v-2$ in $C$ is twice of some part in $D$, then we keep it in $C$.
Otherwise, we split the part $4v-2$ into two equal parts and put
them into $D$. Repeating the procedure until any part of the form $4v-2$ is
 twice of some part in $D$. Denote the triple of partitions we have obtained
  this way by $(A|C'|D')$. It can be checked that such a triple  satisfies the following conditions:

(1) $a_1=c'_1=c'_2=\cdots=c'_{i-1}=c'_i+2=\cdots=c'_{\frac{k-3}{2}}+2$;

(2) $c'_i-c'_{i+\frac{k-3}{2}}\geq 4$, and $2$ appears at most $\frac{a-3}{2}$ times in $C'$;

(3) If there is a part $4v-2$ in $C'$, then the part $2v-1$ must appear in $D'$
either once or twice. If  $4v-2$ does not appear in $C'$ but both $4v$ and $4v-4$
appear in $C'$, then $2v-1$ must appear in $D'$.

\noindent Step 4. For the above conditions (1), (2) and (3), we shall apply Gordon's involution to the pairs $(A|C')$  with $k$ and $a$ replaced by $\frac{k-1}{2}$ and $\frac{a-1}{2}$ respectively and with every part doubled.
The two types of fixed points are as follows
\begin{equation}\label{t5}(A|C'|D')=(4n,4n-2,\ldots,2n+2|(4n)^{\frac{k-a}{2}}\ (4n-2)^{\frac{a-3}{2}}\ (4n-4)^{\frac{k-a}{2}}\ldots|D'),\end{equation}
and
\begin{equation}\label{t6}(A|C'|D')=(4n-2,\ldots,2n|(4n-2)^{\frac{a-3}{2}}\ (4n-4)^{\frac{k-a}{2}}\ (4n-6)^{\frac{a-3}{2}}\ldots|D').\end{equation}

\noindent Step 5. To compute the generating function for the above two types of fixed points,
we need to do some transformations. Set the initial value of $E$ to
 $\emptyset$.  Recall that
 all the parts in $D'$ are odd. We shall distribute the parts of $D'$ to $E$ and $C'$.  To be precise, if a part $d'_i$ is bigger than $2n$, then we move $d_i'$
to $E$. So the remaining  parts in $D'$ are less than $2n$. If there
are two equal parts are left in $D'$, we move one of them to
$E$. Eventually, $D'$ becomes the partition $(2n-1,2n-3,\ldots,3,1)$.
Now we continue to split every part of $C'$ into two equal parts and  move the remaining parts in  $D'$ to $C'$.
It can be checked that the partition $E$ can be any  partition with distinct odd parts. So the fixed points in (\ref{t5}) and (\ref{t6}) take the following equivalent forms
\begin{equation}
(4n,4n-2,\ldots,2n+2|(2n)^{k-a}\ (2n-1)^{a-2}\ (2n-2)^{k-a}\ldots|E),\end{equation}
and
\begin{equation}(4n-2,\ldots,2n|(2n-1)^{a-2}\ (2n-2)^{k-a}\ (2n-3)^{a-2}\ldots|E).\end{equation}
 It is now easy to see that the generating function for the fixed points equals
\[(-q;q^2)_{\infty}
\sum_{n=-\infty}^{\infty}(-1)^nq^{(k+1)n^2+(k+1-a)n}.
\]
This completes the proof. \qed

For example, let $N=189$,
  $k=a=9$ and let $(A|B)=(16,14,12,10\ |\ 9,7^8,5^8,3^8,1^8)$.
  The first three steps are as follows:
\begin{align*}&(16,14,12,10\ |\ 9,7^8,5^8,3^8,1^8)\Rightarrow(16\ |\ 14,12,10\ |\ 14^4,10^4,6^4,2^4\ |\ 9)\\
\Rightarrow&(16,14,12,10\ |\ 14^3,10^3,6^3,2^3\ |\ 9,7,7,5,5,3,3,1,1).
\end{align*}
One can check that the triple of partitions $(16,14,12,10\ |\ 14^3,10^3,6^3,2^3\ |\ 9,7,7,5,5,3,3,1,1)$ is a fixed point of type (\ref{t5}) for $n=4$. Applying the
 transformation in Step 5, we get a triple of partitions
\[(16,14,12,10\ |\ 7^7,5^7,3^7,1^7\ |\ 9,7,5,3,1).\]

\section{Combinatorial proof of  Theorem \ref{thm1.5} }

In this section, we only give a brief description of an involution as
 a combinatorial proof of Theorem \ref{thm1.5} since the idea is similar to
 the involutions for the other two identities of Andrews.

\noindent
{\it Proof of the Theorem \ref{thm1.5}}.
In view of  the relation \[\frac{1}{(-q;q^2)_{\infty}(q;q)_{\infty}}
=\frac{(-q^2;q^2)_{\infty}}{(q^2;q^2)_{\infty}}\]
and Jacobi's triple product identity,
Theorem \ref{thm1.5} can be restated as follows
\begin{equation}\label{c}
(q^2;q^2)_{\infty}\sum_{n=0}^{\infty}\overline{W}_{k,a}(n)q^n=(-q^2;q^2)_{\infty}
\sum_{n=-\infty}^{\infty}(-1)^nq^{(k+1)n^2+(k+1-a)n}.\end{equation}

Let $(A|B)$ be a partition
pair corresponding to the left hand side of (\ref{c}), that is, $A=(
a_1,a_2,\ldots)$ is a
signed partition with distinct even parts and  $B=( b_1,b_2,\ldots)$ is a partition enumerated by $\overline{W}_{k,a}$. We denote the set of such partition  pairs by $\mathcal{Q}_{k,a}$.  We now give  an involution on $\mathcal{Q}_{k,a}$ which
can be described by the following five steps.

\noindent Step 1. We combine two equal parts in $B$ to form an even part. Repeating this
procedure to get an even partition $C=(c_1, c_2, \ldots)$. Suppose that
 the remaining parts of $B$ form a partition $D=(d_1, d_2, \ldots)$.  Since each odd part in $B$ appears an even number of  times, $D$ is a partition with distinct even parts.
We denote  the set of such triples  $(A|C|D)$  by $\mathcal{Q}_{k,a}'$.

\noindent Step 2. Now we compare the largest parts in $A$ and $C$.
 If $c_1> a_1$,  we move the part $c_1$ to
$A$. This leads to  a triple of partitions in $\mathcal{Q}_{k,a}'$.
 If $a_1\geq c_1$, we  move $a_1$ to the partition $C$ to form a triple  in  $\mathcal{Q}_{k,a}'$ unless the following two conditions are satisfied:

(1) $a_1=c_1=c_2=\cdots=c_{i-1}=c_i+2=\cdots=c_{\frac{k-1}{2}}+2$,

(2) $a_1=c_1=c_2=\cdots=c_{i-1}=c_i+2=\cdots=c_{\frac{k-3}{2}}+2$,
and there is an even part $d_j$ in $D$ such that $d_j=a_1/2$ or $ (a_1-2)/2$.

\noindent Step 3. We  wish to construct an involution for further cancelation of
  triples  $(A|C|D)$
in $\mathcal{Q}_{k,a}'$ that satisfy the above
exceptional conditions.
 If
a part $4v$ in $C$ is twice of some part in $D$, then we keep it in $C$.
Otherwise,  we split this part $4v$ into two equal parts and put
them in $D$. After this operation, we  denote the resulting triple of partitions by $(A|C'|D')$.
It can be verified that such a triple $(A|C'|D')$  satisfies the following conditions:

(1) $a_1=c'_1=c'_2=\cdots=c'_{i-1}=c'_i+2=\cdots=c'_{\frac{k-3}{2}}+2$;

(2) $c'_i-c'_{i+\frac{k-3}{2}}\geq 4$, and $2$ appears at most $\frac{a-2}{2}$ times in $C'$;

(3) If there is a part $4v$ in $C'$, then the part $2v$ must appear in $D'$
either once or twice. If  $4v$ does not appear in $C'$ but both $4v-2$ and $4v+2$
appear in $C'$, then $2v$ must appear in $D'$.

\noindent Step 4. For a partition triple $(A|C'|D')$ satisfying  the above conditions (1), (2) and (3), we apply Gordon's involution on the pair
$(A|C')$ with $k$ and $a$ replaced by $\frac{k-1}{2}$ and $\frac{a}{2}$ respectively and with every part doubled. The two types of fixed points are as follows
\begin{equation}\label{t7}(A|C'|D')=(4n,4n-2,\ldots,2n+2|(4n)^{\frac{k-a-1}{2}}
(4n-2)^{\frac{a-2}{2}},(4n-4)^{\frac{k-a-1}{2}}\ldots|D'),\end{equation}
and
\begin{equation}\label{t8}(A|C'|D')=(4n-2,\ldots,2n|
(4n-2)^{\frac{a-2}{2}},(4n-4)^{\frac{k-a-1}{2}}
(4n-6)^{\frac{a-2}{2}}\ldots|D').\end{equation}

\noindent Step 5. We shall make some transformations on the above fixed points
 in order to compute their generating function. Set the initial value of $E$ to $\emptyset$. Since
 all the parts in $D'$ are even, we may distribute the parts of
    $D'$ to $E$ and $C'$. For a fixed point of type \eqref{t7}, if a part $d'_i$ is bigger than $2n+1$, then we move $d_i'$
to $E$. After the completion of this process, the remaining  parts in $D'$ are less than $2n+1$ in the fixed point of type \eqref{t7}. And for a fixed point of type \eqref{t8}, if a part $d'_i$ is bigger than $2n-1$, then we move $d_i'$
to $E$. Then the remaining  parts in $D'$ are less than $2n-1$ in the fixed point of type \eqref{t8}.

If there
are two equal parts that are left in $D'$, we move one of them to
$E$.  Then there are some  remaining parts in $D'$.   For a  fixed point  of  type (\ref{t7}),  the remaining parts in $D'$ are
$(2n,2n-2,\ldots,4,2)$. While for a fixed point  of type (\ref{t8}),
 the remaining parts in $D'$ are  $(2n-2,2n-4,\ldots,4,2)$.

Now we continue to split every part of $C'$ into two equal parts and  move the remaining parts in  $D'$ to $C'$.
It can be seen that  the partition $E$ can be an arbitrary partition
 with distinct odd parts. Therefore, the fixed points can be represented
 by the following types of partitions
\begin{equation}
(4n,4n-2,\ldots,2n+2|(2n)^{k-a}\ (2n-1)^{a-2}\ (2n-2)^{k-a}\ldots|E),\end{equation}
and
\begin{equation}(4n-2,\ldots,2n|(2n-1)^{a-2}\ (2n-2)^{k-a}\ (2n-3)^{a-2}\ldots|E).\end{equation}
Hence we get the generating function
\[(-q^2;q^2)_{\infty}
\sum_{n=-\infty}^{\infty}(-1)^nq^{(k+1)n^2+(k+1-a)n}.
\]
This completes the proof. \qed

For example, let $k=7$, $a=6$ and $N=40$ and
let $(A|B)=(10,2\ |\ 4,4,4,4,4,2,2,2,1,1)$.
 The first  step makes the transformation
\[(10,2\ |\ 4,4,4,4,4,2,2,2,1,1)\Rightarrow(10,2\ |\ 8,8,4,2\ |\ 4,2),\]
One can check that the partition $(10,2\ |\ 8,8,4,2\ |\ 4,2)$ belongs  $U_1^1$ as defined in Section 2 when assuming that every part for Gordon's involution is doubled. Using the map $\beta$, also with every part doubled,  we get  $(10\ |\ 10,8,4,2\ |\ 4,2)$ from the partition pair $(10\ |\ 5,5,4,4,4,2,2,2,1,1)$ via the following steps:
\[(10\ |\ 5,5,4,4,4,2,2,2,1,1)\Rightarrow(10\ |\ 10,8,4,2\ |\ 4,2).\]

\vspace{0.5cm}
 \noindent{\bf Acknowledgments.}  This work was supported by  the 973
Project, the PCSIRT Project of the Ministry of Education,  and the National Science
Foundation of China.

\end{document}